\documentclass[12pt,a4paper]{article}

\usepackage{euler,amsmath,amsfonts,graphicx,url,color,hyperref,defns}
\allowdisplaybreaks

\newcommand{\ibc}{\textsc{ibc}s}
\newcommand{\Res}{\operatorname{Res}}
\newcommand{\KS}{Ku\-ra\-mo\-to--Siva\-shin\-sky}

\newcommand{\bbE}{\mathbb E}
\newcommand{\uij}{u_{i,j}}
\newcommand{\mud}[1][]{\mu_{#1}\delta_{#1}}
\newcommand{\cL}{\mathcal L}

\title{Holistic discretisation ensures fidelity to dynamics in two spatial dimensions}

\author{Tony MacKenzie\thanks{Department of Mathematics and Computing, University of Southern Queensland, Toowoomba, Queensland~4352, Australia.} 
\and 
A.~J. Roberts\thanks{Corresponding author: School of Mathematical Sciences, University of Adelaide, South Australia~5005, Australia. 
\protect\url{mailto:anthony.roberts@adelaide.edu.au}}}

\begin{document}

\maketitle

\begin{abstract}
Developments in dynamical systems theory provides new support for the discretisation of \pde{}s and other microscale systems.  By systematically resolving subgrid microscale dynamics the new approach constructs asymptotically accurate, macroscale closures of discrete models of the \pde.  Here we explore reaction-diffusion problems in two spatial dimensions.  Centre manifold theory ensures that slow manifold, holistic, discretisations exists, are quickly attractive, and are systematically approximated.  Special coupling of the finite elements ensures that the resultant discretisations are consistent with the \pde\ to as high an order as desired.  Computer algebra handles the enormous algebraic details as seen in the specific application to the Ginzburg--Landau equation.  However,  higher order models in 2D appear to require a mixed numerical and algebraic approach that is also developed.  Being driven by the residuals of the equations, the modelling  here may be straightforwardly adapted to a wide class of reaction-diffusion differential and lattice equations in multiple space dimensions.
\end{abstract}

\tableofcontents

\section{Introduction}

Here we extend the dynamical systems `holistic' approach to the macroscale discrete modelling of the class of two dimensional, homogeneous, nonlinear reaction-diffusion equations 
\begin{equation}
\D tu=\nabla\cdot\big[f(u,\nabla u)\nabla u\big]+\alpha g(u)\,.
\label{Erde}
\end{equation} 
Our approach  systematically models subgrid scale processes with the aim of providing an effective closure for the macroscale discretisation.  Sections \ref{S_2D_divide}~and~\ref{sec:nccec} discuss two distinct avenues of theoretical support: that of centre manifold theory, and consistency, respectively.  Such macroscale closures should enable a relatively coarse numerical grid to significantly improve computational speed and stability.

As a particular example, Sections \ref{S_2D_low}~and~\ref{chapnumcm} explore in some detail the real valued, two dimensional, Ginzburg--Landau equation obtained from the \pde~\eqref{Erde} with reaction $g=u-u^3$ and with constant diffusion $f=1$\,, namely
\begin{equation}
\D tu=\nabla^2u+\alpha(u-u^3)\,. \label{E_gl2d}
\end{equation}
We choose this 2D~real Ginzburg--Landau equation as a prototype reaction-diffusion \pde\ because it is well studied and its dynamics well understood~\cite[e.g.]{Gibbon93, Levermore96}.  The non-trivial stable and
unstable steady states of the 2D~Ginzburg--Landau equation~\eqref{E_gl2d}  measure the accuracy and effectiveness of various accuracy models in this application.

The macroscale discretisation is based upon dividing the domain into overlapping finite elements.  Following analogues in one dimension~\cite[e.g.]{Roberts00a,MacKenzie05a}, neighbouring elements are coupled with the non-local
conditions~\eqref{EbcsdL}--\eqref{EsbciL} with the strength of the coupling parametrised by~$\gamma$.   Such coupling of overlapping elements appear analogous to other multiscale methods~\cite[e.g.]{E04, Samaey03b, Gander98}.  Section~\ref{S_2Dcm} then discusses how centre manifold theory assures us of the existence of a slow manifold that is an exactly closed discrete model. Further this slow manifold discretisation is exponentially quickly attractive. Although we cannot find this  exact slow manifold, theory asserts it may be approximated to any asymptotic order in the strength of the interelement coupling~$\gamma$ and the nonlinearity~$\alpha$.  The overlapping finite elements together with the special coupling conditions~\eqref{EbcsdL}--\eqref{EsbciL} assure us that the resultant macroscale discrete models are \emph{also} consistent with the dynamics of the reaction-diffusion \pde\ (Section~\ref{sec:nccec}).

Section~\ref{S_2D_low} outlines the construction, consistency and predictive accuracy of second order asymptotic approximations to the macroscale discretisation of the Ginzburg--Landau \pde~\eqref{E_gl2d}.  To extract another order of accuracy from the algebra, we find (for the first time) the adjoint operator of the diffusion operator on the elements with the nonlocal coupling conditions.  The null space of this adjoint, strikingly similar to a Galerkin basis, enables us to use an integral solvability condition to construct the third order discrete model~\eqref{eq:2dgl3}.

However, higher order holistic models cannot be found analytically.  This inability to construct algebraic approximations is one major difference between systems in one and multiple spatial dimensions.  Section~\ref{chapnumcm} explores how to numerically construct the subgrid field and its evolution in 2D reaction-diffusion equations using the 2D~Ginzburg--Landau equation as an example.  We find that even a relatively coarse subgrid microscale resolution is adequate to accurately predict the macroscale dynamics.

Because of the faithful resolution of subgrid structures and interelement interactions, the resulting discrete models are algebraically complicated as seen, for example, in the third order model~~\eqref{eq:2dgl3}.  Thus users may prefer, as they often do now, to use models of the nonlinear dynamics of lower order.  Then higher order discretisations derived via this approach provide \emph{good} local estimates of the local error in a lower order simulation as it is computed on the fly. 

A further application should be to the `equation-free patch' methodology for simulating multiscale systems~\cite[e.g.]{Kevrekidis03b, Samaey03b}.
This work suggests, analogous to one dimensional dynamics~\cite{Roberts04d, Roberts06d}, that straightforward coupling conditions will also empower efficient and accurate patch dynamics in multiple space dimensions.

\section{Divide the domain into square elements}
\label{S_2D_divide}

\begin{figure}
\begin{center}
\small \setlength{\unitlength}{0.26em}
\begin{picture}(80,80)
\thicklines 
\multiput(10,0)(20,0){4}{
  \multiput(0,0.5)(0,3){27}{\line(0,1){2}}}
\multiput(0,10)(0,20){4}{
  \multiput(0.5,0)(3,0){27}{\line(1,0){2}}}
\multiput(20,20)(20,0){3}{
\multiput(0,0)(0,20){3}{\circle*{2}} } 
\put(12,15){$i-1,j-1$}
\put(34,15){$i,j-1$} \put(52,15){$i+1,j-1$} \put(14,35){$i-1,j$}
\put(38,35){$i,j$} \put(54,35){$i+1,j$} \put(12,55){$i-1,j+1$}
\put(34,55){$i,j+1$} \put(52,55){$i+1,j+1$}
\put(30,6){\vector(1,0){20}}
\put(50,6){\vector(-1,0){20}}
\put(40,3){$h$}
\color{blue}
\multiput(20,20)(40,0){2}{\line(0,1){40}}
\multiput(20,20)(0,40){2}{\line(1,0){40}}
\end{picture}
\end{center}
\caption{The discretisation of a 2D domain into square elements:
The $i,j$th~element (solid) is centred upon the grid point $(x_i,y_j)$; $E_{i,j}$~overlaps neighbouring domains to extend to the neighbouring grid points.}
\label{ch5f2dst}
\end{figure}
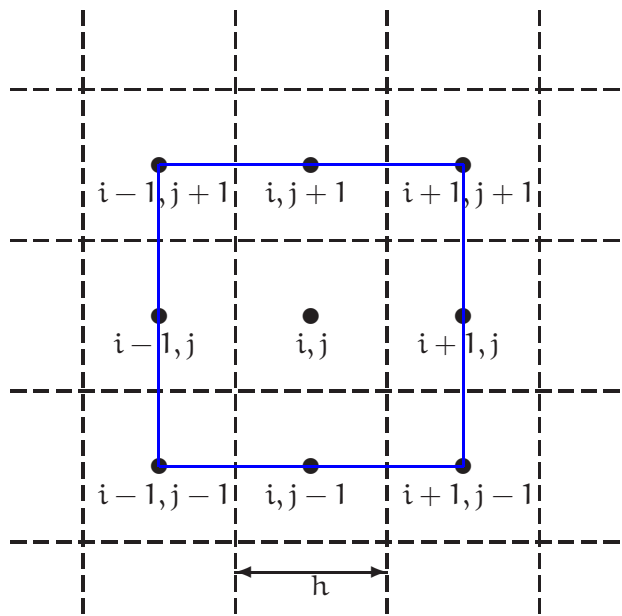

We place the discrete modelling of  two dimensional, reaction-diffusion equations within the purview of centre manifold theory by dividing the domain into overlapping square elements, as shown schematically in Figure~\ref{ch5f2dst}.  The discretisation is similar to that for shear dispersion in a long thin channel~\cite{MacKenzie03}.  The significant difference here is that the discretisation of elements is in both spatial dimensions, not just along the channel as for the shear dispersion application.

\subsection{Extend the non-local IBCs to 2D}

In this initial study, for two dimensional, reaction-diffusion equations we divide the domain into a set of \emph{overlapping} square elements, Figure~\ref{ch5f2dst}.  Define a grid of points~$(x_i,y_j)$ with, for simplicity, constant spacing~$h$. The $i,j$th element,~$E_{i,j}$, is then centred upon~$(x_i,y_j)$ and of width $\Delta x=\Delta y = 2h$\,.  Let $v_{i,j}(x,y,t)$ denote the field in the $i,j$th~element and so evolves according to the reaction-diffusion \pde~\eqref{Erde}; that is, \begin{equation}
\D t{v_{i,j}}=\nabla\cdot\big[f(v_{i,j},\nabla v_{i,j})\nabla v_{i,j}\big]+\alpha g(v_{i,j})\,.
\label{Erdev}
\end{equation}
The original field~$u(x,y,t)$ is then prescribed by $v_{i,j}(x,y,t)$ when $(x,y)\in E_{i,j}$. 

The evolution of the field over the whole domain then depends upon how the elements are coupled together. To couple the dynamics of each element to its neighbours we use coupling `internal boundary conditions' (\ibc) around the $i,j$th~element of
\begin{align}
v_{i,j}(x_{i\pm1},y,t)&=
\gamma v_{i\pm1,j}(x_{i\pm1},y,t)+(1-\gamma)v_{i,j}(x_i,y,t)
,\quad  |y-y_j|<h\,,\quad \label{EbcsdL}\\
v_{i,j}(x,y_{j\pm1},t)&=
\gamma v_{i,j\pm1}(x,y_{j\pm1},t)+(1-\gamma)v_{i,j}(x,y_j,t),
\quad |x-x_i|<h\,.\quad \label{EsbciL}
\end{align}
These \ibc\ are a natural extension to 2D of \ibc\ established for 1D dynamics~\cite[e.g.]{Roberts00a}.   The crucial feature is: with $\gamma=0$ the elements are effectively isolated from each other, dividing the domain into decoupled elements with consequently independent dynamics; whereas with $\gamma=1$ these \ibc\ ensure sufficient continuity between elements to recover the original problem over all space.  Such coupling of overlapping elements appears analogous to the `border regions' of the heterogeneous multiscale method~\cite[e.g.]{E04}, to the `buffers' of the gap-tooth scheme~\cite[e.g.]{Samaey03b}, and to the overlapping domain decomposition that improves convergence in waveform relaxation of parabolic \textsc{pde}s~\cite[e.g.]{Gander98}.

For definiteness in theoretical support, let there be $m$~elements in the domain with the field required to be periodic in both $x$~and~$y$.  For example, the elements may form a $\sqrt m\times\sqrt m$ grid in the domain (any factorisation of~$m$ is feasible).  In principle, the 2D elements could be any shapes, regular or irregular; square elements appear to be easiest to start with.  In this initial work we also avoid issues of \emph{physical domain} boundary conditions on the 2D domain by adopting doubly periodic solutions. (Physical domain boundary conditions have been explored for 1D domains~\cite[e.g.]{Roberts01b, MacKenzie03}.)  These adoptions enable straightforward theoretical statements of support.

\subsection{Centre manifold theory supports discrete models}
\label{S_2Dcm}
This section describes in detail how the \ibc{}~\eqref{EbcsdL}--\eqref{EsbciL} lead to centre manifold theory supporting an accurately closed, discrete model for reaction-diffusion systems~\eqref{Erde}  via its dynamics~\eqref{Erdev} on overlapping elements.

A homotopy in the coupling parameter~$\gamma$ connects the physically relevant discretisation to a tractable base. When parameters $\alpha=\gamma=0$ both the reaction and the coupling on the right-hand side of the \ibc~\eqref{EbcsdL}--\eqref{EsbciL} disappear.  The elements are then effectively isolated from each other and so the resultant diffusion in the \pde~\eqref{Erdev} is particularly simple: exponentially quickly in time, the field~$v_{i,j}$ becomes independently constant within each element.  We use this family of piecewise constant solutions as a basis for analysing the case when the elements are coupled together, $\gamma\neq0$\,.  Particularly interesting is the approximation for full coupling, $\gamma=1$\,, when the \pde~\eqref{Erde} is effectively restored over the whole domain because \ibc~\eqref{EbcsdL}--\eqref{EsbciL} then ensure sufficient continuity between adjacent elements as described previously for 1D \pde{}s~\cite[e.g.]{Roberts98a, MacKenzie05a, Roberts06d}.

The support of centre manifold theory is based upon a linear picture of the dynamics.  Adjoin the dynamically trivial equations
\begin{equation}
    \frac{\partial \alpha}{\partial t}=\frac{\partial \gamma}{\partial t}=0\,,
    \label{EtrivL}
\end{equation}
and consider the reaction-diffusion dynamics in the extended state space $(v_{i,j}(x,y),\gamma,\alpha)$.  In this extended state space, points $\alpha=\gamma=0$ and $v_{i,j}={}$constant are equilibria of the diffusion~\eqref{Erdev}, hence these form a subspace of equilibria, $\bbE_0=\{(v_{i,j},0,0)\}$, in the extended state space.  Linearized about each of the equilibria in~$\bbE_0$, the \pde\ for perturbations~$v'_{i,j}(x,y,t)$ within each element is then the constant coefficient diffusion \pde\ 
\begin{equation}
    \frac{\partial v'_{i,j}}{\partial t}=f_{i,j}\nabla^2v'_{i,j}
    \quad\text{for }(x,y)\in E_{i,j}\text{ for each }i,j,
    \label{eq:ccdpde}
\end{equation}
where the constant diffusivities $f_{i,j}=f(v_{i,j},0)$.  These \pde{}s are decoupled because they are to be solved with the $\gamma=0$ \ibc
\begin{align}&
v'_{i,j}(x_{i\pm1},y,t)=v'_{i,j}(x_i,y,t),\quad |y-y_j|<h\,,
\nonumber\\&
v'_{i,j}(x,y_{j\pm1},t)=v'_{i,j}(x,y_j,t),\quad |x-x_i|<h\,.
    \label{eq:ccdbc}
\end{align}
Thus the following linear eigenmodes are associated with each element:
\begin{align*}&
\alpha=\gamma=0\,,
\\&
v'_{i,j}\propto e^{\lambda_{i,j,k,l} t} \times \cos\big[k\pi(x-X_{i-1/2})/h\big] \times \cos\big[l\pi(y-Y_{j-1/2})/h\big] ,
\end{align*}
inside the $i,j$th~element for all integers $k,l=0,1,2,\ldots$, where the decay rate of each mode is
\begin{equation}
    \lambda_{i,j,k,l}=-f_{i,j}\frac{(k^2+l^2)\pi^2}{h^2}\,;
    \label{Eeigen}
\end{equation}
together with the two trivial modes that firstly $\gamma={}$constant and $\alpha=v'_{i,j}=0$, and that secondly $\alpha={}$constant and $\gamma=v'_{i,j}=0$\,.   In a spatial domain with $m$~elements and when all diffusivities~$f_{i,j}>0$\,, then all eigenvalues are negative, $-f_{i,j}\pi^2/h^2$ or less, except for $m+2$ zero eigenvalues.  Of the $m+2$ zero eigenvalues, one is associated with each of the $m$~elements and two come from the trivial equations~\eqref{EtrivL} for the parameters. That is, the slow subspace is $\{(v_{i,j},\gamma,\alpha)\}$ for constant~$v_{i,j}$. The above arguments establish the following corollary of a centre manifold existence theorem (\cite[p.281]{Carr83b} or~\cite[p.96]{Vanderbauwhede89}).
\begin{corollary}[Existence]
Provided the nonlinear diffusivity~$f$ and reaction~$g$ in~\eqref{Erdev} are sufficiently smooth, and all $f_{i,j}>0$ then  a $m+2$~dimensional slow manifold~${\cal M}$ exists for~\eqref{Erdev}--\eqref{EtrivL} in some finite neighbourhood of the subspace~$\bbE_0$ of equilibria.\footnote{Keep clear the distinction between centre manifold theory and the slow manifolds discussed here: the theory applies to systems where the \emph{real part} of the eigenvalues of critical modes are zero; whereas here we explore and construct slow manifolds because here the eigenvalues are \emph{precisely} zero.}
\end{corollary} 
The slow manifold ${\cal M}$ is parametrized both by the two parameters $\gamma$~and~$\alpha$, and by a measure of the field in each element; we use the grid value $\uij(t)=v_{i,j}(x_i,y_j,t)$ to measure the field in the $i,j$th~element. Using $\vec u$ to denote the vector of such parameters, we write the slow manifold~${\cal M}$  as
\begin{equation}
    v_{i,j}=v_{i,j}(x,y;\vec u,\gamma,\alpha).
    \label{EcmvL}
\end{equation}
These functions~$v_{i,j}(x,y;\vec u,\gamma,\alpha)$, that Sections \ref{S_2D_low}~and~\ref{chapnumcm} construct for the Ginzburg--Landau \pde~\eqref{E_gl2d}, resolve the subgrid scale physical structures as a function of the neighbouring grid values in~$\vec u$.   On this slow manifold~${\cal M}$ the grid values~$\uij$ evolve deterministically
\begin{equation}
    d \uij/dt=\dot u_{i,j}=g_{i,j}(\vec u,\gamma,\alpha)\,,
    \label{EcmgL}
\end{equation}
where $g_{i,j}$~is the restriction of~\eqref{Erdev}--\eqref{EtrivL} to the slow manifold~${\cal M}$.  In essence, this closure of the grid scale dynamics comes from the accurate resolution of the subgrid scale structures.  

Using the value of the field  at the grid points to parametrise the slow manifold provides the necessary `amplitude conditions' to close the problem:
\begin{equation}
\uij=v(x_i,y_j;\vec u,\gamma,\alpha).
\label{Eampl}
\end{equation}  
Many other amplitude conditions are possible such as defining the `amplitudes'~$\uij$ to be the mean field over the $i,j$th~element. However, using the grid values are simple, traditional, and have a direct physical meaning.

Centre manifold theorems~\cite[e.g.]{Carr83b, Vanderbauwhede89} also support the following crucial emergence and approximation properties.
\begin{corollary}[Emergence and approximation] \label{cor:rc}
Provided the nonlinear diffusivity~$f$ and reaction~$g$
in~\eqref{Erdev} are sufficiently smooth, and all $f_{i,j}>0$\,, then  
\begin{itemize}
    \item every solution of the reaction-diffusion dynamics~\eqref{Erdev}--\eqref{EtrivL} that stays within a neighbourhood of the slow manifold~$\cal M$,~\eqref{EcmvL}, approaches exponentially quickly a solution of the discrete model~\eqref{EcmgL} on the slow manifold~\eqref{EcmvL}; and

    \item the order of error in asymptotically approximating the slow manifold and its evolution, \eqref{EcmvL}--\eqref{EcmgL}, is the same as the order of residuals of the governing equations~\eqref{Erdev}--\eqref{EtrivL}.  In particular, because the base equilibria form a subspace~$\bbE_0$, here the approximation is global in the grid values~$\uij$, it is \emph{local only} in the two parameters $\gamma$~and~$\alpha$.
\end{itemize}
\end{corollary}

\section{A slow manifold discretisation of the Ginzburg--Landau equation}
\label{S_2D_low}
We now explore a slow manifold discrete model for a specific reaction-diffusion system, the 2D~Ginzburg--Landau equation equation~\eqref{E_gl2d}.  Substituting \eqref{EcmvL} and \eqref{EcmgL}, the {\textsc{pde}} we solve to form the model is obtained by equating the \pde\ for $\partial v_{i,j}/\partial t$ to that obtained by the chain rule:
\begin{equation}
\D t{v_{i,j}}=\sum_{k,l}\D{u_{k,l}}{v_{i,j}} g_{k,l}
=\nabla^2v_{i,j}+\alpha \left(v_{i,j}-v_{i,j}^3\right)\,. \label{EcmpdeL}
\end{equation}
To construct the slow manifold~\eqref{EcmvL}--\eqref{EcmgL} by solving the \pde~\eqref{EcmpdeL} with coupling and amplitude conditions involves considerable algebraic detail. These algebraic details of the construction  are handled straightforwardly by iteration in computer algebra~\cite[e.g.]{Roberts96a}. The specific procedure used here, documented elsewhere~\cite{MacKenzie09a}, solves the equations using iteration to drive to zero the residuals of the governing differential equation~\eqref{EcmpdeL} and its interelement coupling \ibc~\eqref{EbcsdL}--\eqref{EsbciL}.  Hence by the Approximation Corollary~\ref{cor:rc} we construct correspondingly accurate approximations to the slow manifold of~\eqref{EcmpdeL}.  These approximations, upon setting coupling parameter $\gamma=1$\,, form 2D discrete models of the the 2D~Ginzburg--Landau \pde~\eqref{E_gl2d}.

One consequence of using computer algebra is that there is no need to record in this article most of the considerable algebraic detail in constructing the models.  Those wishing to verify the correctness of the results recorded herein should download and examine the corresponding technical report~\cite{MacKenzie09a} that details the precise computer algebra procedure.  Because the algorithm is based upon driving the residuals to zero, the critical aspect of the procedure is simply the correct coding of the computation of the residuals of the governing equations.

\paragraph{The $\Ord{\gamma^3+\alpha^3}$ holistic discretisation}
\label{S_2d_hol}
Satisfying the \pde\ and \ibc{} to residuals of~$\Ord{\gamma^3+\alpha^3}$ the computer algebra procedure~\cite[\S2.2]{MacKenzie09a} gives subgrid fields which are too complex to record here. The corresponding evolution of the grid values on the slow manifold are
\begin{align}
\dot u_{i,j}=&\frac{\gamma}{h^2}\delta^2 \uij +\alpha\left(\uij - \uij^3\right)
\nonumber\\&{}
-\frac{\gamma^2}{12h^2}\delta^4 \uij
+\alpha \gamma \left(\rat{1}{12} \delta^2\uij^3-\rat14 \uij^2 \delta^2\uij\right)
+\Ord{\gamma^3+\alpha^3}\,,
\label{E_gl2d_g1a1}
\end{align}
where the centred difference operator applies in both spatial dimensions,
\begin{eqnarray*}
\delta^2\uij&=&u_{i+1,j}+u_{i-1,j}+u_{i,j+1}+u_{i,j-1}-4\uij\,,
\\
\delta^4\uij&=&u_{i+2,j}+u_{i-2,j}+u_{i,j+2}+u_{i,j-2}
\\&&{}
-4(u_{i+1,j}+u_{i-1,j}+u_{i,j+1}+u_{i,j-1}) +12\uij\,,
\end{eqnarray*}
The model~\eqref{E_gl2d_g1a1} is simply the extension to two spatial dimensions of the $\Ord{\gamma^3+\alpha^3}$ holistic model of the 1D Ginzburg--Landau equation~\cite{MacKenzie05}.

The holistic discrete model has the dual justification of consistency with the \pde\ in addition to the justification provided by centre manifold theory.  As proven in Section~\ref{sec:nccec}, consistency for such discrete models follows from the coupling \ibc~~\eqref{EbcsdL}--\eqref{EsbciL}~\cite{Roberts00a}. Set the coupling parameter $\gamma=1$ in the discrete equation~\eqref{E_gl2d_g1a1} to recover the holistic discrete model of the Ginzburg--Landau \pde~\eqref{E_gl2d} in 2D.  To test consistency, we expand the finite differences of~\eqref{E_gl2d_g1a1} in a Taylor series in the grid spacing~$h$~\cite[\S2.5]{MacKenzie09a} to find the equivalent continuum \pde{} for the $\Ord{\gamma^3+\alpha^3}$ holistic model~\eqref{E_gl2d_g1a1} is
\begin{equation}
\D tu=\nabla^2u+\alpha(u-u^3)
+\frac{h^2\alpha}{2}u|\nabla u|^2
- \frac{h^4}{90} \left( \Dn x6u+\Dn y6u\right)
+\Ord{\alpha^3+h^6}\,.
\label{E_gl2d_epde}
\end{equation}
The $\Ord{\gamma^3+\alpha^3}$ holistic model is $\Ord{h^6+\alpha^3}$ consistent,
maintaining in 2D the dual justification of holistic discretisation found for 1D \pde{}s~\cite{Roberts98a, Roberts00a}.  Section~\ref{sec:nccec} proves this consistency in some generality for 2D \pde{}s.

\subsection{Illustration of the subgrid field in 2D}
\label{S_2D_perf}
\begin{figure}
 \centering
\includegraphics[width=\textwidth]{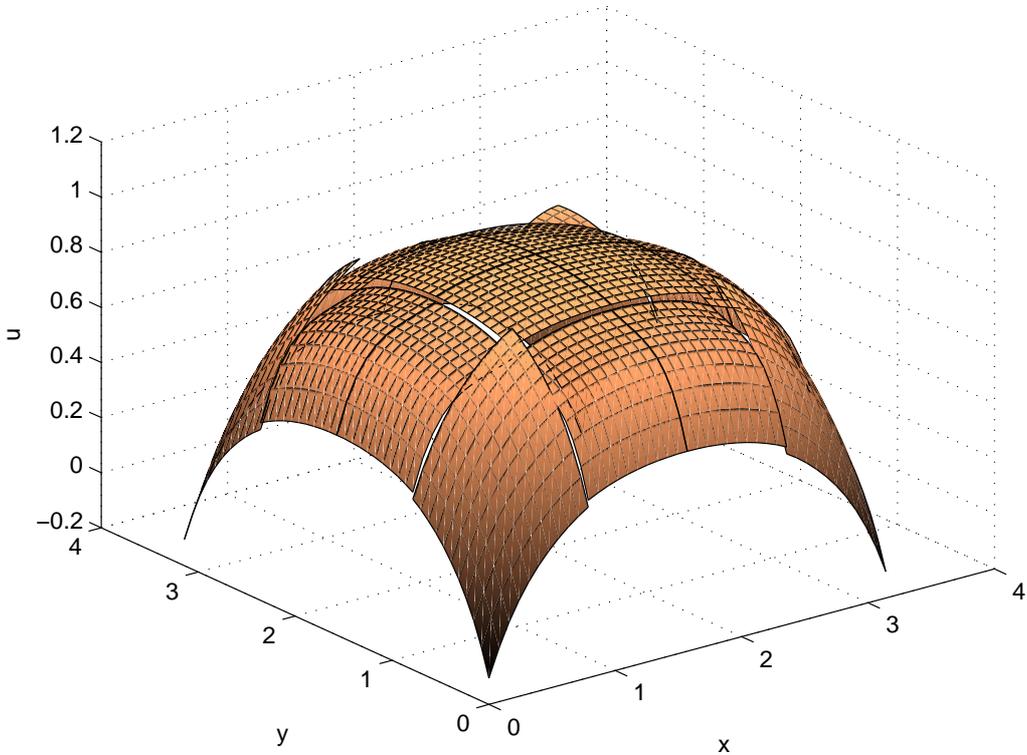}
\caption{an example of the subgrid field for the $\Ord{\gamma^2,\alpha^2}$~holistic
model with $4\times4$ elements on $[0,\pi]\times[0,\pi]$ at nonlinearity $\alpha=6$}
    \label{ch5fig3}
\end{figure}

Here we plot the subgrid fields for a coarse grid solution of the $\Ord{\gamma^2,\alpha^2}$~holistic model (obtained from~\eqref{E_gl2d_g1a1} by omitting the $\gamma^2$~term and then evaluating at $\gamma=1$) of the 2D~Ginzburg--Landau equation~\eqref{E_gl2d}.  We restrict attention to a doubly odd symmetric solution that is $2\pi$--doubly periodic.  That is,
\begin{align}&
u(x,y,t)=u(x+2\pi,y,t), \quad u(x,y,t)=-u(2\pi-x,y,t)
\nonumber\\&
u(x,y,t)=u(x,y+2\pi,t), \quad u(x,y,t)=-u(x,2\pi-y,t)\,.
\label{E_gl_bc_po_2d}
\end{align}
Figure~\ref{ch5fig3} shows the subgrid fields for the $\Ord{\gamma^2,\alpha^2}$~holistic model with $4\times4$ elements on $[0,\pi]\times[0,\pi]$ at nonlinearity $\alpha=6$.  The subgrid fields exhibit the nonlinear subgrid structure of the 2D~Ginzburg--Landau equation and its interaction through the \ibc{}.  The subgrid fields are comprised of actual solutions, albeit approximate, of the 2D~Ginzburg--Landau \pde.

Note the subgrid fields have noticeable jumps across the boundaries of the elements.  Higher order holistic models should reduce these jumps across the boundaries as seen for the holistic models of the \KS{} equation~\cite{MacKenzie05a}.

\paragraph{The $\Ord{\gamma^2,\alpha^2}$ holistic model needs improving}

We investigate the performance of the $\Ord{\gamma^2,\alpha^2}$ holistic model on coarse grids by comparing its bifurcation diagram to an accurate solution.  Again we restrict our attention to doubly odd symmetric solutions that are $2\pi$-doubly periodic~\eqref{E_gl_bc_po_2d}.

The bifurcation information is calculated using the continuation software \textsc{auto}~\cite{auto} and \textsc{xppaut}~\cite{xppaut} as outlined  for the \KS{} equation in MacKenzie's PhD dissertation~\cite{MacKenzie05}.  In such bifurcation diagrams the blue curves indicate stable steady state solutions and red curves indicate unstable steady state solutions.  The open squares indicate steady state bifurcations.

\begin{figure}
\centering
\includegraphics[width=.9\textwidth]{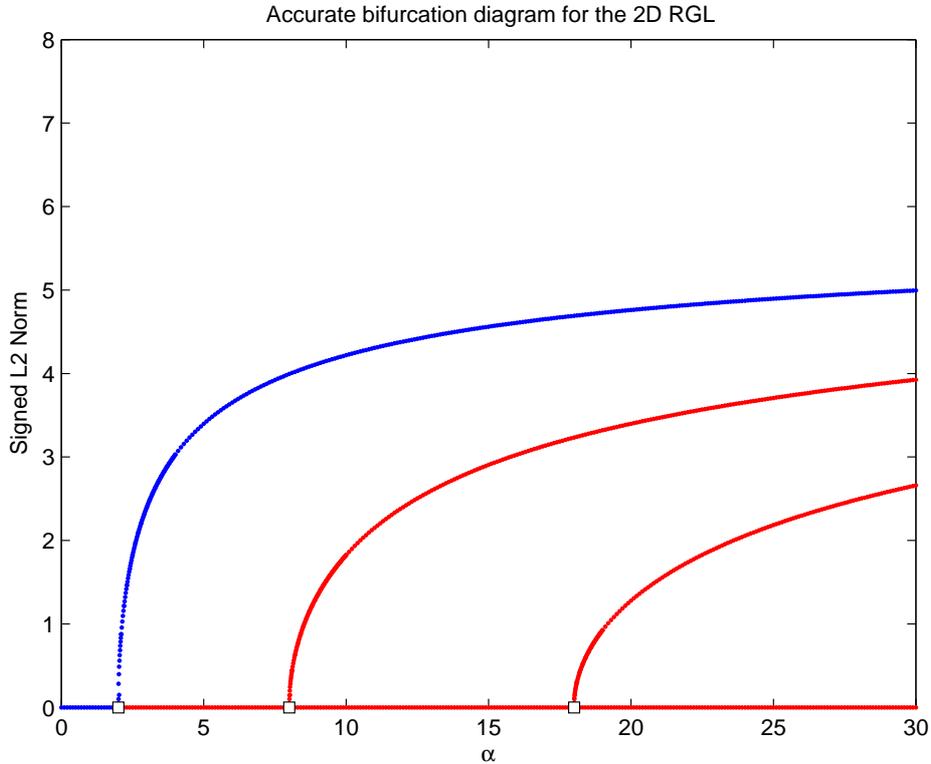}
\caption{Accurate bifurcation diagram of the 2D~Ginzburg--Landau equation system for $0\le\alpha\le30$.  It is constructed with a fourth order centered difference approximation with $24\times24$ points on $[0,\pi]\times[0,\pi]$.  Blue curves indicate stable steady state solutions and red curves indicate unstable steady state solutions.  The open squares indicate steady state bifurcations.} \label{ginz2d_acc}
\end{figure}

Figure~\ref{ginz2d_acc} shows an accurate bifurcation diagram of the 2D~Ginzburg--Landau \pde.  It is constructed with a fourth order centered difference approximation with $24\times24$ points on $[0,\pi]\times[0,\pi]$.  The trivial solution undergoes steady state bifurcations at $\alpha=2,8,18$ leading to the unimodal, bimodal and trimodal branches respectively.  For $1<\alpha<30$, only the unimodal branch is stable and all other branches are unstable.

\begin{figure}
\centering
\includegraphics[width=\textwidth]{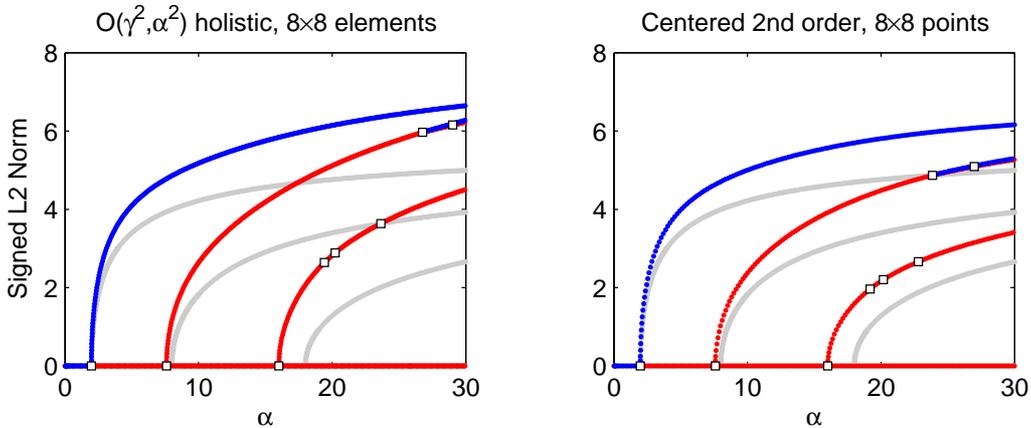}
\caption{Bifurcation diagrams for the 2D~Ginzburg--Landau equation with $8\times8$ elements on $[0,\pi]\times[0,\pi]$ for (a)~${\cal O}(\gamma^2,\alpha^2)$ holistic model (b)~second order centered difference.  The accurate bifurcation diagram is shown in grey.}
\label{ginz2d_g1r1_8_bif}
\end{figure}

Figure~\ref{ginz2d_g1r1_8_bif} shows a comparison of the ${\cal O}(\gamma^2,\alpha^2)$ holistic model and a second order explicit centered difference approximation with $8\times8$ elements on $[0,\pi]\times[0,\pi]$.  The accurate bifurcation diagram is shown in grey.
The ${\cal O}(\gamma^2,\alpha^2)$ holistic model does not perform as well as the second order centered difference approximation; this is similar to performance observed for the one dimensional Ginzburg--Landau \pde~\cite{MacKenzie05}.

\paragraph{Higher order models need numerical construction}
To improve the accuracy of the holistic discretisation we need to compute higher orders in either coupling~$\gamma$, or nonlinearity~$\alpha$, or both.   Improved accuracy occurs at higher order in comparable 1D problems~\cite{MacKenzie05a}.    However, apparently it is not possible to \emph{analytically} construct higher order subgrid fields in 2D: apparently the subgrid fields required for our closures are no longer in the class of multivariate polynomials.  Instead, numerical methods must be used to find the subgrid fields as described in Section~\ref{chapnumcm}.  However, the well known `solvability condition' in asymptotic mathematical methods empowers us to derive the next order in the evolution, analytically from the residuals, without needing to find the next order of the subgrid fields.

\subsection{The adjoint provides an extra order of accuracy}
\label{sec:apeoa}

We scrounge an extra order of accuracy from the `solvability condition'~\cite[e.g.]{Pavliotis07} applied to residuals of the next asymptotic order.  Because the linear operator used to find corrections to the subgrid field is singular---the operator necessarily has homogeneous solutions that compose the slow subspace~$\mathbb E_0$---the Fredholm alternative is that the `right-hand side' of the equation for the subgrid fields must lie in the range of the singular operator.  This solvability condition is enough to determine an extra correction to the evolution.  

Recall from linear algebra that to be in the range of the operator, the solvability condition is that the right-hand side must be orthogonal to the null space of the adjoint operator.  Thus the first task of this section is to find the adjoint operator of the linear constant diffusion \pde~\eqref{eq:ccdpde} with its boundary conditions~\eqref{eq:ccdbc}.  Second, we find a basis for the null space.  Lastly, more computer algebra provides the required extra order in the evolution.

\begin{figure}
\centering
\setlength{\unitlength}{1ex}
\begin{picture}(27,25)(0,0)
\thicklines
\put(-1,-2){
\multiput(5,5)(10,0){3}{\line(0,1){20}}
\multiput(5,5)(0,10){3}{\line(1,0){20}}
\put(5.5,10){\put(0,0){$L$}\put(10,0){$M$}\put(20,0){$R$}}
\put(10,5.5){\put(0,0){$B$}\put(0,10){$C$}\put(0,20){$T$}}
\put(4.5,3){\put(0,0){$X_{i-1}$}\put(10,0){$X_i$}\put(20,0){$X_{i+1}$}}
\put(0.5,4.5){\put(0,0){$Y_{j-1}$}\put(0,10){$Y_j$}\put(0,20){$Y_{j+1}$}}
}
\end{picture}
\caption{each element is effectively divided into four subregions by the non-locality of the boundary conditions~\eqref{eq:ccdbc}.  To derive the adjoint, label the edges of these four subregions as shown.}
\label{fig:radjoint}
\end{figure}
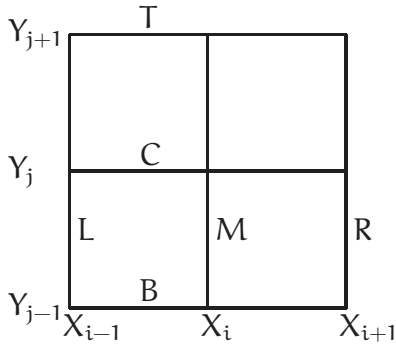

The decoupling of the elements, provided by $\gamma=0$ in the boundary conditions~\eqref{eq:ccdbc}, simplifies finding the adjoint: we need only consider each element in isolation.  Thus define the inner product to be the integral over the $i,j$th~element:
\begin{displaymath}
\left<v,w\right>=\iint_{E_{i,j}}vw\,dx\,dy
=\int_{Y_{j-1}}^{Y_{j+1}} \int_{X_{i-1}}^{X_{i+1}} vw\,dx\,dy \,.
\end{displaymath}
To find the adjoint recognise that each element is subdivided into four subregions, shown schematically in Figure~\ref{fig:radjoint}, by the coupling of the boundary values to internal values by the boundary conditions~\eqref{eq:ccdbc}.  In addition, there exists previously implicit conditions that the subgrid field~$v$ and its gradient are continuous throughout the element. Then integration by parts, or the divergence theorem, transforms the inner product
\begin{align*}
\left<\delsq v,w\right>&
=\left<v,\delsq w\right>
-\int_L v_xw-vw_x\,dy  +\int_R v_xw-vw_x\,dy
\\&\qquad{} 
+\int_{M-} v_xw-vw_x\,dy -\int_{M+} v_xw-vw_x\,dy
\\&\qquad{}+\text{analogous integrals on $B$, $C$~and~$T$,}
\end{align*}
where specific parts of the boundary integrals are labelled as shown in Figure~\ref{fig:radjoint}. Using superscripts to denote evaluation, continuity requires $v^{M\pm}=v^M$ and $v_x^{M\pm}=v_x^M$\,, and the boundary conditions~\eqref{eq:ccdbc} imply $v^L=v^M=v^R$\,.  Thus the inner product
\begin{align*}
\left<\delsq v,w\right>
&=\left<v,\delsq w\right>
+\int_{Y_{j-1}}^{Y_{j+1}} \Big[-v_x^Lw^L+v^Mw_x^L  +v_x^Rw^R-v^Mw_x^R
\\&\qquad{} 
+v_x^Mw^{M-} -v^Mw_x^{M-} -v_x^Mw^{M+} +v^Mw_x^{M+}
\Big] dy
\\&\qquad{}+\text{(analogous $x$~integrals of $y$~derivatives)}
\\&=\left<v,\delsq w\right>
+\int_{Y_{j-1}}^{Y_{j+1}} \Big[-v_x^Lw^L +v_x^Rw^R 
+v_x^M\left(w^{M-}-w^{M+}\right)
\\&\qquad{} 
+v^M\left(w_x^L-w_x^R -w_x^{M-} +w_x^{M+}\right)
\Big] dy
\\&\qquad{}+\text{(analogous $x$~integrals of $y$~derivatives).}
\end{align*}
For the adjoint, these integrals on the right-hand side must vanish for all smooth fields~$v$. Consequently, the null space of the adjoint operator satisfies Laplace's equation~$\delsq w=0$ with conditions: firstly, that $w$~is zero around the edges~$L$, $R$, $B$~and~$T$ of the element; secondly, that $w$~is continuous on the interior partitions $M$~and~$C$ (but its gradients may be discontinuous there); thirdly, that $w_x^L -w_x^{M-} +w_x^{M+}-w_x^R=0$\,; and lastly, that $w_y^B -w_y^{C-} +w_y^{C+}-w_y^T=0$\,.

\begin{figure}
\centering
\setlength{\unitlength}{1mm}
\begin{picture}(110,78)(0,0)
\put(2,0){\includegraphics{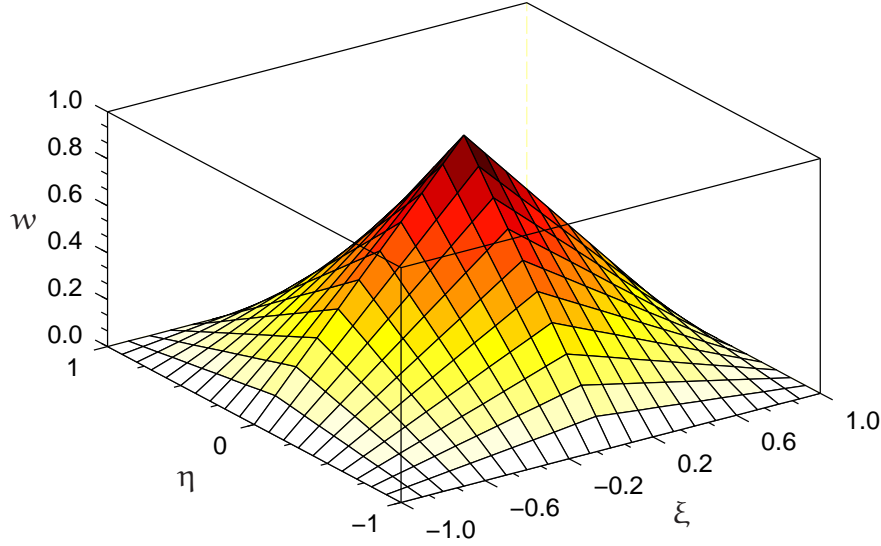}}
\put(22,10){$\eta$}
\put(87,5){$\xi$}
\put(0,44){$w$}
\end{picture}
\caption{basis `pyramid'  for the null space of the adjoint operator on an element: $w=(1-|\xi|)(1-|\eta|)=(1-|x-X_i|/h)(1-|y-Y_i|/h)$.}
\label{fig:zadjoint}
\end{figure}

Because of these conditions, the null space of the adjoint is spanned by the `pyramid' $w=(1-|x-X_i|/h)(1-|y-Y_i|/h)$ as displayed in Figure~\ref{fig:zadjoint}. The solvability condition is then that the integral of the subgrid residuals of the \pde\ with this~$w$ over the $i,j$th~element determines a correction to the model evolution.  

It will not escape your notice that the solvability condition integral parallels integrals in the Galerkin finite element method.  Thus the Galerkin finite element method may be viewed as a leading approximation to our systematic slow manifold closure of discrete modelling.

\paragraph{Return to the discrete Ginzburg--Landau model}
Computer algebra readily computes the subgrid residuals of the Ginzburg--Landau \pde~\eqref{E_gl2d} to the next higher order (detailed elsewhere~\cite[\S2.3]{MacKenzie09a}).  Taking the inner product with the adjoint null vector~$w$, and remembering contributions from the inter-element coupling conditions~\eqref{EbcsdL}--\eqref{EsbciL}, gives the discrete model, in gory detail,
\begin{align}
\dot u_{i,j}=&\frac{\gamma}{h^2}\delta^2 \uij 
-\frac{\gamma^2}{12h^2}\delta^4 \uij
+\frac{\gamma^3}{90h^2}\delta^6 \uij
\nonumber\\&{}
+\alpha\left(\uij - \uij^3\right)
+ \frac{\gamma \alpha}{12} \left( \delta^2\uij^3-3 \uij^2 \delta^2\uij\right)
\nonumber\\&{}
+\frac{\gamma^2\alpha}{720}\left[
\uij^2( 222\delta^2\uij +24\delta^4\uij -3\delta_x^2\delta_y^2\uij )
\right. \nonumber\\&\quad\left.{}
+\uij( -102\delta^2\uij^2 
+36\{\delta^2\uij\}^2 
+6\{\delta_x^2\uij\}\{\delta_y^2\uij\} 
-144\{\mud\uij\}^2  )
\right. \nonumber\\&\quad\left.{}
-6\{\mud[y]\delta_x^2\uij\}\{\mud[y]\uij^2\}
-6\{\mud[x]\delta_y^2\uij\}\{\mud[x]\uij^2\}
\right. \nonumber\\&\quad\left.{}
+12\{\mud\uij^2\}\{\mud\uij\}
+12\{\mud^3\uij\}\{\mud\uij^2\}
-\rat32\{\delta^2\uij^2\}\{\delta_x^2\delta_y^2\uij\}
\right. \nonumber\\&\quad\left.{}
+3\{\delta^4\uij\}\{\delta^2\uij^2\}
-3\{\delta_x^2\uij^2\}\{\delta_y^2\uij\}
-3\{\delta_y^2\uij^2\}\{\delta_x^2\uij\}
\right. \nonumber\\&\quad\left.{}
+9\{\delta^2\uij^2\}\{\delta^2\uij\}
-8\delta^4\uij^3 -6\delta^2\uij^3 +\delta_x^2\delta_y^2\uij^3
\right]
\nonumber\\&{}
+\frac{\gamma\alpha^2h^2}{240}\left[ 3\uij^4\delta^2\uij
+6\uij^2\delta^2(\uij-\uij^3) 
-2\delta^2\uij^3 +3\delta^2\uij^5 \right]
\nonumber\\&{}
+\Ord{\gamma^4+\alpha^4}\,,
\label{eq:2dgl3}
\end{align}
This complicated macroscale discrete closure arises from resolving subgrid scale nonlinear dynamics within the finite elements.

To independently check the above model, consider the small element limit as $h\to0$\,.   Upon setting the coupling parameter $\gamma=1$ we expect to recover classic consistency.   Straightforward computer algebra~\cite[\S2.5]{MacKenzie09a} finds that for small grid size~$h$ the equivalent partial differential equation to the discrete model~\eqref{eq:2dgl3} is
\begin{align}
\D tu=&\delsq u+\alpha(u-u^3)
+\frac{h^6}{560}\left[\partial_x^8u+\partial_y^8u\right]
\nonumber\\&{}
+\frac{\alpha h^4}{60}\left[
uu_{xy}^2 +2u_xu_yu_{xy}
-8(u_x^2u_{xx}+u_y^2u_{yy})
\right.\nonumber\\&\left.\quad{}
-5u(u_{xx}^2+u_{yy}^2)
-14u(u_xu_{xxx}+u_yu_{yyy})
\right]
\nonumber\\&{}
+\Ord{\alpha^4+h^8}
\label{eq:eqpde48}
\end{align}
As found for 1D problems with analogous element coupling~\cite{Roberts00a}, this slow manifold discrete model is consistent to the 2D Ginzburg--Landau equation~\eqref{E_gl2d} to high order in the element size~$h$, both for the linear diffusion and the nonlinear reaction.  Section~\ref{sec:nccec} proves our coupling conditions construct consistent discrete models in general.  But remember that Corollary~\ref{cor:rc} independently provides strong support for the relevance of the model~\eqref{eq:2dgl3} at the \emph{finite} element sizes used in simulations.

However, the discrete model~\eqref{eq:2dgl3} is as high an order of accuracy as we can construct analytically.  The next Section~\ref{chapnumcm} shows how to numerically solve for the subgrid scale field in order to construct the macroscale discrete model. It serves as a proof of principle for applying the holistic method to \pde{}s of two or more spatial dimensions.

\section{Generally compute 2D subgrid fields numerically}
\label{chapnumcm}

The holistic discretisation of \pde{}s is based upon centre manifold theory~\cite[e.g.]{Carr81,Carr83b, Kuznetsov95} to resolve the subgrid fields and hence more accurately close the macroscale discrete model. Here the subgrid field is constructed numerically for the Ginzburg--Landau equation~\eqref{E_gl2d} in 2D.  

New complexities arise.  Although the spatial structure is obtained numerically, the slow manifold, subgrid field is also parametrised by the grid values~$u_{i,j}$, the interelement coupling parameter~$\gamma$, and the nonlinear parameter~$\alpha$.  Therefore,  the construction involves symbolic parameters.  The algorithm required to develop the holistic model must efficiently solve the corresponding mixed numerical and symbolic problem. \emph{The focus of this section is on this novel numerical construction of the subgrid fields and not the performance of the holistic models.} 

Numerical construction of the subgrid field introduces errors which are separate from the orders of errors of the holistic model.  These errors from the numerical construction of the holistic discretisation are the major concern. The numerical construction of the subgrid field and its evolution has challenging details: \S\ref{S_num_2D_ext},~How should the subgrid problem be solved?  \S\ref{S_num_2D_perf},~What subgrid resolutions will accurately reproduce the analytical holistic models?  \S\ref{S_num_eff},~What is an efficient implementation?

\subsection{Outline the numerical slow manifold in 2D}
\label{S_num_2D_ext}
\begin{figure}
\begin{center}
\small \setlength{\unitlength}{0.26em}
\begin{picture}(80,80)
\thicklines 
\multiput(10,0)(20,0){4}{
  \multiput(0,0.5)(0,3){27}{\line(0,1){2}}}
\multiput(0,10)(0,20){4}{
  \multiput(0.5,0)(3,0){27}{\line(1,0){2}}}
\multiput(20,20)(20,0){3}{
\multiput(0,0)(0,20){3}{\circle*{2}} } 
\put(12,15){$i-1,j-1$}
\put(34,15){$i,j-1$} \put(52,15){$i+1,j-1$} \put(15,35){$i-1,j$}
\put(38,35){$i,j$} \put(55,35){$i+1,j$} \put(12,55){$i-1,j+1$}
\put(34,55){$i,j+1$} \put(52,55){$i+1,j+1$}
\put(30,6){\vector(1,0){20}}
\put(50,6){\vector(-1,0){20}}
\put(40,3){$h$}
\color{blue}
\multiput(25,20)(5,0){7}{\line(0,1){40}}
\multiput(20,25)(0,5){7}{\line(1,0){40}}
\multiput(20,25)(40,0){2}{\line(0,1){30}}
\multiput(25,20)(0,40){2}{\line(1,0){30}}
\end{picture}
\end{center}
\caption{Example of the $4 \times 4$ interval subgrid in 2D; note that such a subgrid labelled as ``$4\times4$'' actually extends to be $9\times9$ when overlapped with neighbouring elements as shown.}
\label{fig:2D_sub}
\end{figure}
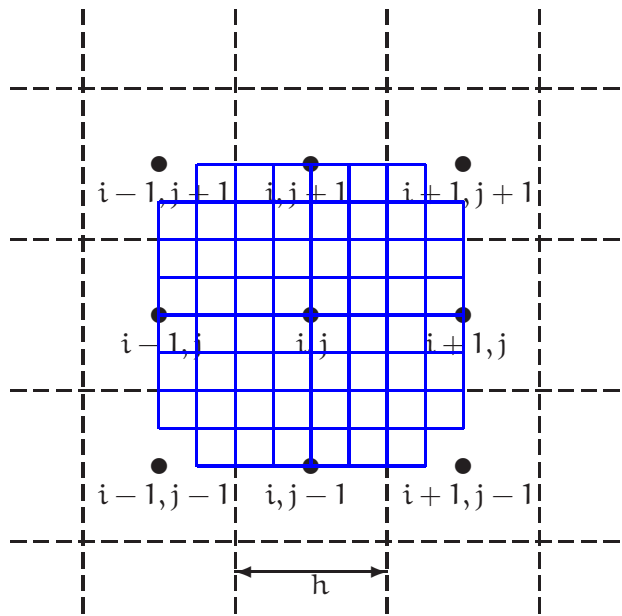

In each element we discretise the microscale subgrid as shown in Figure~\ref{fig:2D_sub}.  At each subgrid grid point we seek the evolution of $u_{i,j,k,\ell}$ for subgrid indices $|k|,|\ell|<n$ where, for example, Figure~\ref{fig:2D_sub} shows $n=4$\,. The subgrid is shown solid (blue) for this particular example of a $4\times4$ interval subgrid.  The subgrid field extends to $X_{i\pm1}$~and~$Y_{j\pm1}$ in order to allow the application of the 2D non-local \ibc{} \eqref{EbcsdL}--\eqref{EsbciL}.  The subgrid field does not extend to the extreme corners because the subgrid discrete Laplacian applied to the interior points does not involve the subgrid field at the corners.
 
After discretising the subgrid field in the $i,j$th~element, classic finite differences approximate the spatial derivatives of the subgrid field of the Ginzburg--Landau \pde~\eqref{E_gl2d}:
\begin{align}
\dot v_{i,j,k,\ell}
=&\frac{v_{i,j,k+1,\ell}+v_{i,j,k-1,\ell}+v_{i,j,k,\ell+1}+v_{i,j,k,\ell-1}
-4v_{i,j,k,\ell}}{(h/n)^2} 
\nonumber\\&{}
+\alpha(v_{i,j,k,\ell}-v_{i,j,k,\ell}^3).
\label{eq:nde}
\end{align}
These microscale discretised equations are solved at each of the subgrid points inside each of the elements. The elements are coupled with \ibc\ analogous to~\eqref{EbcsdL}--\eqref{EsbciL}, namely
\begin{align}
&v_{i,j,\pm n,\ell}=\gamma v_{i\pm1,j,0,\ell}+(1-\gamma)v_{i,j,0,\ell}
\quad |\ell|<n\,,
\label{eq:nibcx}\\
&v_{i,j,k,\pm n}=\gamma v_{i,j\pm1,k,0}+(1-\gamma)v_{i,j,k,0}
\quad |k|<n\,.
\label{eq:nibcy}
\end{align}
The same centre manifold theorems apply to the system~\eqref{eq:nde}--\eqref{eq:nibcy} to assure us of the  existence, relevance and construction of a slow manifold, macroscale discrete model of the dynamics. The macroscale slow manifold to construct is that the subgrid field $v_{i,j,k,\ell}=v_{i,j,k,\ell}(\vec u,\gamma,\alpha)$ where, defining $u_{i,j}=v_{i,j,0,0}$\,, the macroscale grid values~$\vec u$ evolve according to $\dot u_{i,j}=g_{i,j}(\vec u,\gamma,\alpha)$.

We employ an iteration scheme to find the microscale subgrid field and the macroscale slow evolution.  The initial approximation is that of a constant field in each element: $v_{i,j,k,\ell}\approx u_{i,j}$ such that $\dot u_{i,j}=g_{i,j}\approx 0$\,.  Given any current approximation, $\vec v_{i,j}$~and~$g_{i,j}$, we seek an improved approximation $\vec v_{i,j}:=\vec{v}_{i,j}+\vec  v_{i,j}'$ and $g_{i,j}:={g}_{i,j}+g_{i,j}'$ where $\vec v_{i,j}'$~and~$g_{i,j}'$ are corrections to be found in each iteration.    At each iteration, the following linear equations driven by the current residuals are solved for the corrections~$\vec{v}_{i,j}'$ and~$g_{i,j}'$
\begin{align}
\left(\frac nh\right)^2\left(v'_{i,j,k+1,\ell}+v'_{i,j,k-1,\ell}
+v'_{i,j,k,\ell+1}+v'_{i,j,k,\ell-1}
-4v'_{i,j,k,\ell}\right) -g_{i,j}'&=\Res_{\ref{eq:nde}}\,,\nonumber\\
v'_{i,j,\pm n,\ell}-v'_{i,j,0,\ell}&=\Res_{\ref{eq:nibcx}}\,,\nonumber\\
v'_{i,j,k,\pm n}-v'_{i,j,k,0}&=\Res_{\ref{eq:nibcy}}\,,\nonumber\\
v'_{i,j,0,0}&=0\,.
\label{E_glnum_eqs}
\end{align}
The iteration repeats until all residuals are zero to a specified order of error. This iteration scheme follows that for the analytic construction of the subgrid field and is documented in full detail in a separate technical report~\cite[\S3]{MacKenzie09a}. Centre manifold theory then assures us that the resultant macroscale discrete model $\dot u_{i,j}=g_{i,j}(\vec u,\gamma,\alpha)$ is accurate to the same order of error.

For example, for the coarsest possible subgrid field, with just $n=2$ subintervals, the numerical discretisation of the subgrid field gives a low order, macroscale model as
\begin{align}
\dot u_{i,j}=&\frac{\gamma}{h^2}\delta^2 \uij +\alpha\left(\uij - \uij^3\right)
\nonumber\\&{}
-\frac{\gamma^2}{16h^2}\delta^4 \uij
+\alpha \gamma \left(\rat1{16} \delta^2\uij^3-\rat3{16} \uij^2 \delta^2\uij\right)
+\Ord{\gamma^3+\alpha^3}\,,
\label{E_num_gl_g1r1_2g}
\end{align}
Compare with the analytic macroscale discrete model~\eqref{E_gl2d_g1a1}: the terms in the first line are identical; the same higher order terms appear in the second line but the coefficients are in error by~$25$\%.  This correspondence is promising for such a coarse microscale subgrid discretisation.

\subsection{Low resolution subgrids are accurate in 2D}
\label{S_num_2D_perf}

How do macroscale models constructed via a numerical microscale, such as~\eqref{E_num_gl_g1r1_2g}, compare with analytic macroscale models? We compare in two ways: one via the convergence of the coefficients; and the other by the accuracy of the predicted bifurcation diagrams.  It appears that the microscale subgrid need not be of high resolution.

\begin{table}
\caption{coefficients in the $\Ord{\gamma^3+\alpha^3}$ models, such as~\eqref{E_num_gl_g1r1_2g}, evidently converge to the correct analytic coefficients, in~\eqref{E_gl2d_g1a1} and labelled~$\infty$ in the table, with errors~$\Ord{1/n^2}$ as the resolution of the microscale grid improves.}
\label{tbl:num_gl_g1r1_2g}
\begin{center}
\begin{tabular}{cccc}
$n$ & 
$\gamma^2\delta^4\uij/h^2$ & 
$\alpha\gamma\delta^2\uij^3$ & 
$\alpha\gamma\uij^2\delta^2\uij$ \\ \hline
\vphantom{$\Big|$}
$2$ & $-\frac1{16}$ & $\frac1{16}$ & -$\frac3{16}$ \\
\vphantom{$\Big|$}
$4$ & $-\frac5{64}$ & $\frac5{64}$ & -$\frac{15}{64}$ \\
\vphantom{$\Big|$}
$8$ & $-\frac{21}{256}$ & $\frac{21}{256}$ & -$\frac{63}{256}$ \\
\hline
\vphantom{$\Big|$}
$\infty$ & $-\frac1{12}$ & $\frac1{12}$ & -$\frac1{4}$ \\
\hline
\end{tabular}
\end{center}
\end{table}

First look at the coefficients of the~$\Ord{\gamma^3+\alpha^3}$ models such as~\eqref{E_num_gl_g1r1_2g}.  Recall that the number of microscale subgrid points, from one macro-grid point to the next, in each dimension, is~$n$.  The coefficients linear in the coupling parameter~$\gamma$ and nonlinearity~$\alpha$ are exact for all $n\geq 2$\,, only the higher order coefficients vary with subgrid resolution.  Table~\ref{tbl:num_gl_g1r1_2g} tabulates coefficients in these nonlinear terms, those of~$\Ord{\gamma^2+\alpha^2}$ in models such as~\eqref{E_num_gl_g1r1_2g}, for some values of~$n$.  Evidently the coefficients converge to the exact values with error~$\Ord{1/n^2}$.  We expect such quadratic convergence from the quadratic  modelling in~\eqref{eq:nde} of the subgrid scale dynamics.

\begin{table}
\caption{maximum errors in the coefficients of the $\Ord{\gamma^4+\alpha^4}$ model~\eqref{eq:2dgl3} when approximated numerically at three different subgrid resolutions.  The decrease by at least a factor of four, upon doubling~$n$, indicates quadratic convergence.}
\label{tbl:me2dgl3}
\begin{center}
\begin{tabular}{llllll}
$n$ & $\gamma^2/h^2$ & $\gamma\alpha$ & $\gamma^3/h^2$ & $\gamma^2\alpha$ & $\gamma\alpha^2h^2$ \\ \hline
$2$ & $0.021 $ & $0.062 $ & $0.0033 $ & $0.14 $ & $0.0016   $ \\
$4$ & $0.0052$ & $0.016 $ & $0.00086$ & $0.040$ & $0.000098 $ \\
$8$ & $0.0013$ & $0.0039$ & $0.00022$ & $0.010$ & $0.0000061$ \\ \hline
\end{tabular}
\end{center}
\end{table}

Similarly, we compare numerically obtained coefficients for the $\Ord{\gamma^4+\alpha^4}$ model~\eqref{eq:2dgl3}.  Because of the complexity of the model we make a limited comparison: for each order in $\alpha$~and~$\gamma$ in~\eqref{eq:2dgl3}, Table~\ref{tbl:me2dgl3} reports the largest error in the numerically obtained coefficient for three different subgrid scale resolutions.  Evidently, these maximum errors decrease like~$1/n^2$ to confirm the accuracy of the numerical description of the subgrid scale dynamics. 

Incidentally, this agreement between the numerical model, obtained by resolving the subgrid scale dynamics at order~$\Ord{\gamma^3+\alpha^3}$, and the analytic model~\eqref{eq:2dgl3} supports the derivation in Section~\ref{sec:apeoa} of an extra order of accuracy in the macroscale model without necessarily having to resolve the subgrid scale dynamics at the same order.

\begin{figure}
\centering
\includegraphics[width=\textwidth]{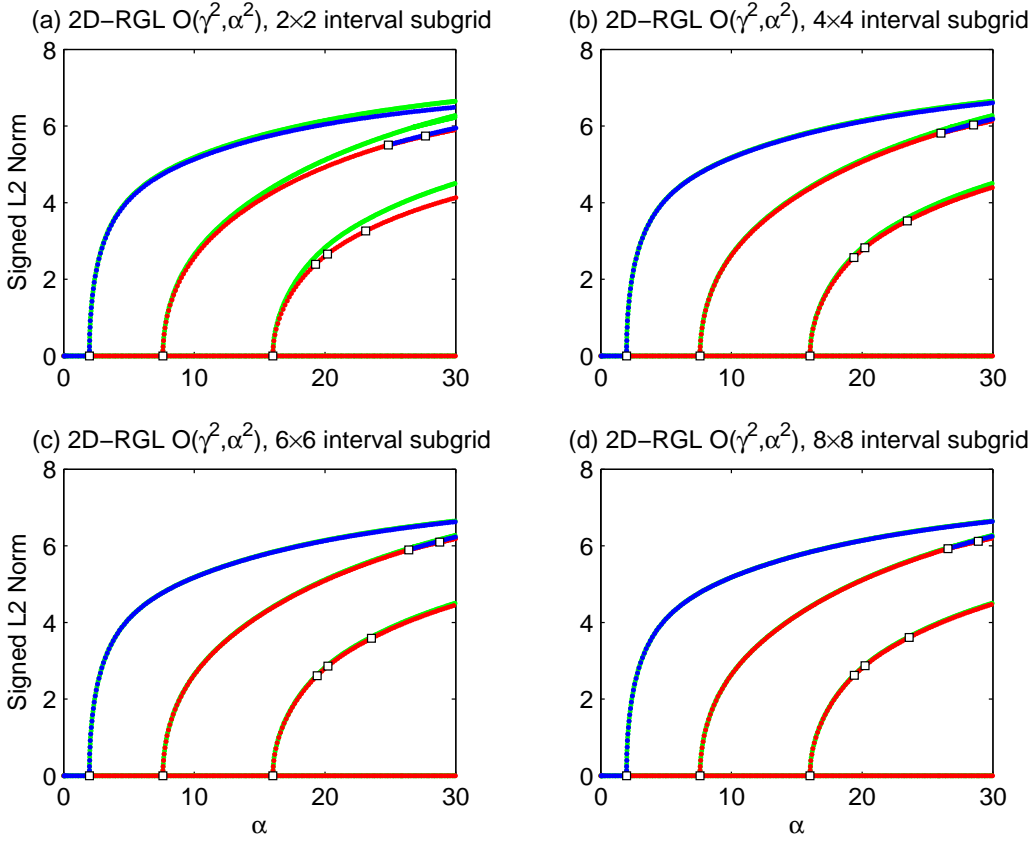}
\caption{Bifurcation diagrams of the $O(\gamma^2,\alpha^2)$ holistic
models of the 2D Ginzburg--Landau system with $8\times 8$ macroscale
elements on $[0,\pi]\times[0,\pi]$ for subgrid resolutions of
(a)~$2\times2$, (b)~$4\times4$, (c)~$6\times6$ and (d)~$8\times8$
intervals.  The bifurcation diagram for the analytically constructed
model is shown in green.}
\label{ginz2d_num_g1r1}
\end{figure}

Second, we turn to the bifurcation diagram to see the sort of errors incurred in using the approximate models. Figure~\ref{ginz2d_num_g1r1} shows the bifurcation diagrams for the $\Ord{\gamma^2,\alpha^2}$ holistic model of the Ginzburg--Landau system for four subgrid resolutions.  Here the equilibria shown in green are not the accurate solution of the Ginzburg--Landau system, but rather the equilibria of the analytic $\Ord{\gamma^2,\alpha^2}$ holistic model in 2D (obtained from~\eqref{E_gl2d_g1a1} by omitting the $\gamma^2$~term).  Observe that with a subgrid resolution of just $4\times4$ intervals the bifurcation diagram for the numerically constructed $\Ord{\gamma^2,\alpha^2}$ holistic model is almost indiscernible from the analytic model over nonlinearity $0\le\alpha\le20\,$.  Higher subgrid resolutions are indiscernible to even larger nonlinearity~$\alpha$.

\begin{figure}
\centering
\includegraphics[width=\textwidth]{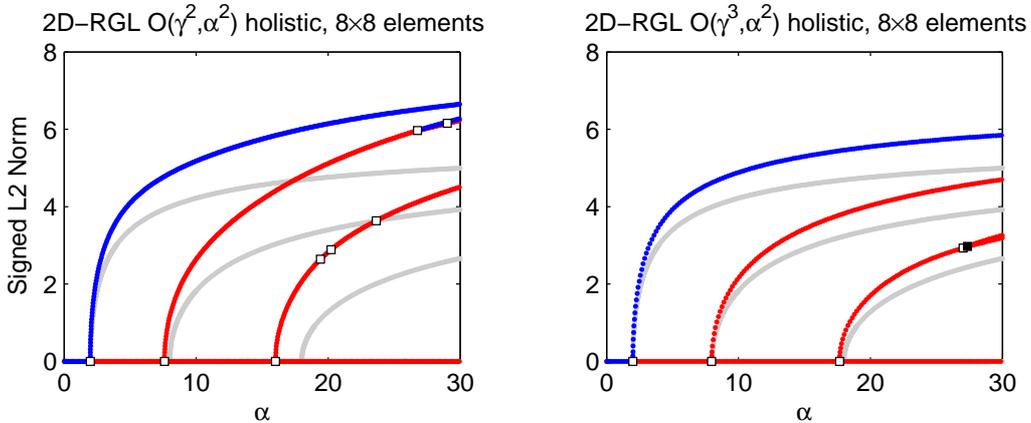}
\caption{Bifurcation diagrams of the (a)~$\Ord{\gamma^2,\alpha^2}$ and (b)~$\Ord{\gamma^3,\alpha^2}$ holistic models of the 2D Ginzburg--Landau system with $8\times8$ elements on $[0,\pi]\times[0,\pi]$.  The accurate bifurcation diagram without stability information is shown in grey.}
\label{ginz2d_num_8_bif}
\end{figure}

As a last comparison of bifurcation diagrams, Figure~\ref{ginz2d_num_8_bif} shows one example confirming  that by computing to higher order in the interelement interactions, and resolving the subgrid scale structures numerically, the predictions of the numerically derived models do improve. 

Numerically resolving the microscale subgrid structures does generate usefully accurate, slow manifold, macroscale discretisations.

\subsection{An efficient computer algebra approach is crucial}
\label{S_num_eff}

The difficulty associated with the numerical construction of the subgrid field is the mixed discrete and symbolic nature of the equations involved in the iteration scheme.  The size of the system of equations increases as the subgrid resolution improves, and the complexity of the symbolic nonlinear residuals increases quickly as higher order holistic models are constructed.  

Computer algebra packages such as \textsc{reduce}~\cite{Reduce04} or \textsc{mathematica}~\cite{Wolfram96} have general routines that solve systems of equations such as~\eqref{E_glnum_eqs}.  However, these \verb|solve| routines were inefficient for the many equations and complicated expressions involved with better subgrid resolution and higher order holistic models.  Even a low resolution, $n=2$ interval subgrid, took many minutes in both \textsc{reduce} and \textsc{mathematica} using their built in \verb|solve| routines.  Instead we develop an approach that is practical for implementation with a large number of complicated symbolic terms.

\paragraph{Transform to constant coefficient}
Recall that at each step of the iteration scheme we solve a problem for updates to the subgrid field $\vec{v}'_{i,j}$ and its evolution~$g_{i,j}'$.  Multiply the first (field) equation in~\eqref{E_glnum_eqs} by~$h^2$ and replace~$g'_{i,j}$ by $\mathcal G'=h^2g'_{i,j}$\,.  Then the left-hand side of the new form of the equations has numerical constant coefficients; algebraic expressions only occur in the right-hand side.

Further, the left-hand side of the new equations remain the same for every iteration.  Consequently, the first iteration constructs an $LU$~factorisation of the left-hand side, which is then used to solve equation~\eqref{E_glnum_eqs} for updates in every iteration~\cite[\S3]{MacKenzie09a}. The $LU$~decomposition is performed once and requires approximately $\frac13 N^3$~operations~\cite[e.g.]{Press92}.  Here the number of equations for the subgrid structure are $N=(2n+1)^2+1$\,; for example, $N=290$ for $n=8$ microscale intervals in the subgrid. At each step of the iteration scheme the \textsc{lu} factorisation algorithm operates on the symbolic residual vector.  Perhaps 2D and 3D problems could be solved more efficiently through iterative multigrid~\cite{mccormick92,brigg00} or incomplete $LU$~factorisation and Krylov subspace methods~\cite{kelley95,vorst03}. Such alternatives remain for later exploration.

\begin{table}
\caption{\textsc{reduce} and \textsc{mathematica} computational times for numerical construction of $\Ord{\gamma^4,\alpha^2}$ holistic models of the one dimensional Ginzburg--Landau equation for various subgrid scale resolutions,~$n$.}
\label{table_comptime_red}
\begin{center}
\begin{tabular}{crr}
$n$&\textsc{reduce}&\textsc{mathematica}\\
\hline
$2$&$1.1$\,s&$70.2$\,s\\
$4 $&$3.1$\,s&$215.4$\,s\\
$8 $&$8.3$\,s&$367.6$\,s\\
$16 $&$23.7$\,s&$517.7$\,s\\
\hline
\end{tabular}
\end{center}
\end{table}

\paragraph{Reduce was faster}
Computational comparison experiments found that the computer algebra \textsc{reduce} was an order of magnitude faster than \textsc{mathematica}. Table~\ref{table_comptime_red} lists the computational time for the \textsc{reduce} and the \textsc{mathematica} implementation for constructing $\Ord{\gamma^4,\alpha^2}$~holistic models of the \emph{one dimensional} Ginzburg--Landau equation with subgrid resolutions of~$2$, $4$, $8$ and~$16$ subgrid intervals.  These times were observed on a Pentium~III, 750\,MHz processor, with 256\,Mb~\textsc{ram}, running \textsc{reduce}~3.7, under Windows~XP. Table~\ref{table_comptime_red} shows the \textsc{reduce} implementation was $20$--$70$ times faster than the \textsc{mathematica} implementation  (despite the repeated help of the \textsc{mathematica} news group).  Thus we used \textsc{Reduce}.

\textsc{reduce} is an interpreted language and as such it interprets each command and allocate memory as each command is executed at runtime.  It would be possible to write a purpose built compiled program in some mainstream language to study any specific \pde.  This option is not considered in the scope of this article.

\section{Non-local coupling conditions enforce consistency}
\label{sec:nccec}

Recall that the constructed holistic models of the Ginzburg--Landau dynamics are consistent with the \pde\  as the grid size $h\to0$\,, see equations \eqref{E_gl2d_epde}~and\eqref{eq:eqpde48} in Section~\ref{S_2D_low}.   Now we prove that general consistency follows from the specific choice of nonlocal interelement coupling conditions~\eqref{EbcsdL}--\eqref{EsbciL}.  

We start with a similar theorem to one previously proved for the consistency of holistic discretisation in one space dimension~\cite{Roberts00a}.  The critical difference here is in the proof: previously the proof was constructive whereas here it is not.  Avoiding a constructive proof has two consequences: it is essential here as we do not know analytic forms for the slow manifold subgrid fields in 2D; and the new proof easily caters for nonlinear reaction.  Because the theorem here is more powerful, an immediate corollary proves consistency of a 2D holistic discretisation with the 2D linear \pde.

\begin{theorem}[1D consistency] \label{thm:1dc}
Consider the \pde\ $\partial_tu=\cL u +g(u)$ for some local, isotropic, homogeneous, linear operator~$\cL$, and for some smooth nonlinear reaction~$g$. Model the dynamics on overlapping elements of an equi-spaced grid $X_i=ih$\,.  Let $v_i(x,t)$~denote the subgrid field in the $i$th~element satisfying the \pde\ $\partial_tv_i=\cL v_i +g(v_i)$ on the interval~$(X_{i-1},X_{i+1})$ with the moderated interelement coupling conditions
\begin{equation}
v_i(X_{i\pm1},t)=\gamma v_{i\pm1}(X_{i\pm1},t)+(1-\gamma)v_i(X_i,t)\,.
\label{eq:1dibc}
\end{equation}
When interelement interactions are truncated to residuals~$\Ord{\gamma^p}$ the grid values $U_i(t)=v_i(X_i,t)$, at full coupling $\gamma=1$\,, evolve consistently with the \pde\ $\partial_tu=\cL u+g(u)$\,. 
\end{theorem}

\begin{proof}
We proceed with some classic operator algebra~\cite[e.g.]{npl61}.  The principle obstacle is to  transform subgrid spatial differences, indicated by subscript~$x$, into macroscale grid differences, indicated by subscript~$i$. Begin with the \pde\ on the $i$th~element: $\partial_t v_i=\cL v_i +g(v_i)$\,.  Because the operator~$\cL$ is isotropic and homogeneous it may be formally expanded in even centred differences as
\begin{displaymath}
\cL  =\sum_{k=0}^\infty \ell_{2k}\delta_x^{2k} =\ell_0+\ell(\delta_x^2)\,,
\end{displaymath}
for some coefficients~$\ell_{2k}$ and corresponding function~$\ell$.
For example, in application to reaction~diffusion \pde{}s, we would write the diffusion operator $\partial_x^2=\big[\frac2h\sinh^{-1}(\rat12\delta_x)\big]^2$.  Then provided the leading coefficient $\ell_2\neq0$\,, the \pde\ $\partial_tv_i=[\ell_0+\ell(\delta_x^2)]v_i+g(v_i)$ is equivalently written $\ell^{-1}\big[(\partial_t-\ell_0)v_i-g(v_i)\big]=\delta_x^2 v_i$ where $\ell^{-1}$~is the inverse of function~$\ell$ (implicitly assumed invertible).    
Evaluate this last form of the \pde\ at $x=X_i$ so that the $v_i$ on the left-hand side becomes simply the grid value~$U_i$ and by the coupling conditions~\eqref{eq:1dibc}\footnote{If the leading coefficient in the expansion of~$\cL$ is $\ell_{2n}\neq0$\,, because the lower order coefficients are zero (or asymptotically small as in the \KS\ \pde), then more coupling conditions like~\eqref{eq:1dibc} couple with the next nearer neighbouring elements.} the centred spatial difference on the right-hand side becomes the centred grid difference~$\gamma\delta_i^2U_i$\,.  This evaluation then gives the evolution $\ell^{-1}\big[(\partial_t-\ell_0)U_i-g(U_i)\big]=\gamma\delta_i^2 U_i$ on the macroscale grid.  

Now reverting the inverse function, this grid evolution is equivalent to
\begin{equation}
\partial_t U_i=\big[\ell_0+\ell(\gamma\delta_i^2)\big]U_i +g(U_i)
=\sum_{k=0}^\infty \gamma^k\ell_{2k}\delta_i^{2k} U_i +g(U_i)\,.
\label{eq:dgis}
\end{equation}
For example, for the diffusion equation
\begin{align*}
\partial_tU_i&=\frac1{h^2}\left[2\sinh^{-1}\big(\rat12\sqrt\gamma\delta_i\big)\right]^2U_i
\\&{}
=\frac1{h^2}\left[\gamma\delta_i^2-\frac{\gamma^2}{12}\delta_i^4
+\frac{\gamma^3}{90}\delta_i^6 -\frac{\gamma^4}{560}\delta_i^8 +\cdots\right]U_i\,.
\end{align*}
Thus a truncation of~\eqref{eq:dgis} to errors~$\Ord{\gamma^p}$ results in a discrete model with stencil width of $2p-1$\,.  But specifically relevant to the theorem is that the equivalent differential equation of this discrete model evaluated at full coupling:  the error in approximating~$\cL$ by the truncated version of~\eqref{eq:dgis} (at full coupling $\gamma=1$) is dominated by the leading neglected term, namely $\ell_{2p}\delta_i^{2p}$.  As the element size $h\to 0$, this error is~$\Ord{\ell_{2p}h^{2p}}$.  For example, for the diffusion operator~$\partial_x^2$, the coefficients $\ell_{2k}=\Ord{1/h^2}$ and so the discrete model is consistent with the diffusion \pde\ to error~$\Ord{h^{2p-2}}$ as grid size $h\to 0$\,.
\end{proof}

For interest, immediately generalise this theorem to a class of nonlinear reaction-dissipation \pde{}s.

\begin{corollary} \label{thm:1dcn}
Consider the \pde\ $\partial_tu=\cL F(u) +g(u)$ for some local, isotropic, homogeneous, linear operator~$\cL$, and for some smooth nonlinear~$F$ and reaction~$g$. Model the dynamics on overlapping elements of an equi-spaced grid $X_i=ih$\,.  Let $v_i(x,t)$~denote the subgrid field in the $i$th~element satisfying the \pde\ $\partial_tv_i=\cL F(v_i) +g(v_i)$ on~$(X_{i-1},X_{i+1})$ with the moderated interelement coupling conditions
\begin{equation}
F\big[v_i(X_{i\pm1},t)\big]
=\gamma F\big[v_{i\pm1}(X_{i\pm1},t)\big]
+(1-\gamma)F\big[v_i(X_i,t)\big]\,.
\label{eq:1dibcn}
\end{equation}
When interelement interactions are truncated to residuals~$\Ord{\gamma^p}$ the grid values $U_i(t)=v_i(X_i,t)$, at full coupling $\gamma=1$\,, evolve consistently with the \pde\ $\partial_tu=\cL F(u)+g(u)$\,. 
\end{corollary}

\begin{proof}
Naturally generalise the proof for Theorem~\ref{thm:1dc}.
\end{proof}

In this theorem and corollary the subgrid microscale operator~$\cL$ need not be differential.  For example, $\cL$~could be a microscopic discretisation as in the numerical construction of the previous Section~\ref{chapnumcm}: for example,
\begin{equation*}
\delta^2_{\text{subgrid}}
=4\sinh^2\left[\rat12h_{\text{subgrid}}\partial_x\right]
=4\sinh^2\left[\frac h{2n}\partial_x\right]
=4\sinh^2\left[\rat1n\sinh^{-1}\rat12\delta_i\right].
\end{equation*}
The proof still holds.  Further, the proof still holds if the subgrid field is not a scalar.  This last observation empowers the following corollary on general consistency in two dimensions that Section~\ref{S_2D_low} observed specifically for the Ginzburg--Landau \pde..

\begin{corollary}[2D consistency] \label{thm:2dc}
Consider a reaction-diffusion \pde\ $\partial_tu=\nabla^2 u+g(u)$ (such as the Ginzburg--Landau equation~\eqref{E_gl2d}) modelled on overlapping elements with subgrid fields $v_{i,j}(x,y,t)$ coupled by conditions~\eqref{EbcsdL}--\eqref{EsbciL}.
When the interactions are truncated to residuals~$\Ord{\gamma^p}$ the grid values $U_{i,j}(t)=v_{i,j}(X_i,Y_j,t)$, at full coupling $\gamma=1$\,, evolve consistently with the \pde\ to error~$\Ord{h^{2p-2}}$. 
\end{corollary}

\begin{proof}
Apply the previous Theorem~\ref{thm:1dc} twice.  First, treat coordinate~$y$ as a parameter so that $\ell_0=\partial_y^2$ and $\ell(\delta_x^2)=\partial_x^2$.  Then by Theorem~\ref{thm:1dc} the semi-discrete `grid' values $v_{i,j}(X_i,y,t)$, for discrete~$i$ and parametrised by continuous~$y$, evolve consistently with the reaction-diffusion \pde.  Second, treat index~$i$ as a parameter and consider the discrete modelling in coordinate~$y$ so that now $\ell_0$ involves operators acting on the $x$-grid and $\ell(\delta_y^2)=\partial_y^2$.  Then by Theorem~\ref{thm:1dc} the 2D grid values $U_{i,j}=v_{i,j}(X_i,Y_j,t)$ evolve consistently with the semi-discrete system generated in the first step, which in turn is consistent with the diffusion \pde.  The observation at the end of the proof for Theorem~\ref{thm:1dc} provides the order of error.
\end{proof}

\section{Conclusion}

We explored novel macroscale discretisation of reaction-diffusion dynamics in two spatial dimensions.  This work generalises considerable earlier work on modelling dynamics in one dimension.  By dividing space into \emph{overlapping} elements, we have shown that the specific interelement coupling conditions~\eqref{EbcsdL}--\eqref{EsbciL} have important theoretical and practical consequences.

Section~\ref{S_2D_divide} discussed how these coupling conditions ensure that centre manifold theory applies to assure us of the existence, relevance and approximation of the slow manifold that is the macroscale discretisation of general reaction-diffusion equations in two dimensions.  Further, Section~\ref{sec:nccec} proved that the resulting discrete models will also be consistent, as the macroscale grid size $h\to0$\,, with the continuum or microscale dynamics.  Thus the holistic discretisations generated in this novel approach have the dual justification of both consistency and centre manifold theory.

This strong theoretical support appears to be straightforwardly generalisable to reaction-diffusion dynamics in three or more dimensions.  The support also appears to be straightforwardly generalisable to higher order \pde{}s, just as the theory supports the discrete modelling of one dimensional higher order \pde{}s such as the \KS\ \pde~\cite{MacKenzie00a, MacKenzie05a}.  An interesting issue for further research is whether there are alternative interelement coupling conditions in two or more dimensions that have additional desirable properties.  For example, can one find interelement coupling that generates the so-called compact discrete approximations to the diffusion Laplace operator? rather than the spatially extended ones generated here.

Sections \ref{S_2D_low}~and~\ref{chapnumcm}, using the example of the real Ginzburg--Landau \pde, explored many of the technical issues necessary to apply the approach.   In contrast to one spatial dimension, we found that a purely algebraic approach can only model to low order of accuracy.  Although such low order models do reasonably accurately predict the dynamics, as seen in bifurcation diagrams, higher order accuracy is desirable.  Consequently we introduced and explored an approach where the microscale subgrid dynamics are described numerically, but with algebraic expressions for coefficients so that we construct an algebraic model for the macroscale discretisation.   The computer algebra package used does make a difference:  we found \textsc{reduce} more than an order of magnitude faster than \textsc{Mathematica}.  Further research may see if iterative techniques based upon the sparse microscale interactions are more effective than the direct $LU$~factorisation employed here.  This work provides a novel approach and theory for the sound and accurate closure of macroscale discretisations.

\paragraph{Acknowledgement}  The Australian Research Council Discovery Project grants DP0774311 and DP0988738 helped support this research.

\bibliographystyle{plain}
\bibliography{more2,ajr,bib}

\begin{thebibliography}{10}

\bibitem{brigg00}
W.L. Brigg, V.E. Hanson, and S.F. Mc{C}ormick.
\newblock {\em A multigrid tutorial}.
\newblock Society for Industrial and Applied Mathematics ({SIAM}),
  Philadelphia, 2000.

\bibitem{Carr81}
J.~Carr.
\newblock {\em Applications of centre manifold theory}, volume~35 of {\em
  Applied Math Sci}.
\newblock Springer-Verlag, 1981.

\bibitem{Carr83b}
J.~Carr and R.G. Muncaster.
\newblock The application of centre manifold theory to amplitude expansions.
  {II}. infinite dimensional problems.
\newblock {\em J. Diff. Eqns}, 50:280--288, 1983.

\bibitem{auto}
E.J. Doedel, R.C. Paffenroth, A.R. Champneys, T.F. Fairgrieve, Yu.A. Kuznetsov,
  B.~Sandstede, and X.~Wang.
\newblock Auto 2000: Continuation and bifurcation software for ordinary
  differential equations (with {HomCont}).
\newblock Technical report, Caltech, 2001.

\bibitem{E04}
Weinan E, Bjorn Engquist, Xiantao Li, Weiqing Ren, and Eric Vanden-Eijnden.
\newblock The heterogeneous multiscale method: A review.
\newblock Technical report,
  \url{http://www.math.princeton.edu/multiscale/review.pdf}, 2004.

\bibitem{xppaut}
B.~Ermentrout.
\newblock {XPPAUT} 5.0 - the differential equations tool.
\newblock Technical report,
  [\url{http://www.math.pitt.edu/~bard/bardware/xpp_doc.pdf}], 2001.

\bibitem{Gander98}
Martin~J. Gander and Andrew~M. Stuart.
\newblock Space-time continuous analysis of waveform relaxation for the heat
  equation.
\newblock {\em SIAM Journal on Scientific Computing}, 19(6):2014--2031, 1998.

\bibitem{Gibbon93}
J.D. Gibbon.
\newblock Weak and strong turbulence in the complex {Ginzburg-Landau} equation.
\newblock In G.R. Sell, C.~Foais, and R.~Temam, editors, {\em Turbulence in
  fluid flows--A dynamical systems approach}, volume~55 of {\em The {IMA}
  volumes in mathematics and its applications}, pages 33--48. Springer-Verlag,
  1993.

\bibitem{Reduce04}
A.C. Hearn.
\newblock [\url{http://www.reduce-algebra.com}].
\newblock 2004.

\bibitem{kelley95}
C.T. Kelley.
\newblock {\em Iterative methods for linear and nonlinear equations}, volume~16
  of {\em Frontiers in applied mathematics}.
\newblock Society for Industrial and Applied Mathematics {SIAM}, 1995.

\bibitem{Kevrekidis03b}
I.~G. Kevrekidis, C.~W. Gear, J.~M. Hyman, P.~G. Kevrekidis, O.~Runborg, and
  K.~Theodoropoulos.
\newblock Equation-free, coarse-grained multiscale computation: enabling
  microscopic simulators to perform system level tasks.
\newblock {\em Comm. Math. Sciences}, 1:715--762, 2003.

\bibitem{Kuznetsov95}
Y.~A. Kuznetsov.
\newblock {\em Elements of applied bifurcation theory}, volume 112 of {\em
  Applied Mathematical Sciences}.
\newblock Springer--Verlag, 1995.

\bibitem{Levermore96}
C.D. Levermore and M.~Oliver.
\newblock The complex {Ginzburg-Landau} equation as a model problem.
\newblock In P.~Deift, C.D. Levermore, and C.E. Wayne, editors, {\em Dynamical
  systems and probabilistic methods in partial differential equations},
  volume~35 of {\em Lectures in Applied Mathematics}, pages 141--190. American
  Mathematical Society, 1996.

\bibitem{MacKenzie00a}
T.~Mackenzie and A.~J. Roberts.
\newblock Holistic finite differences accurately model the dynamics of the
  {Kuramoto--Sivashinsky} equation.
\newblock {\em ANZIAM~J.}, 42(E):C918--C935, 2000.
\newblock \url{http://anziamj.austms.org.au/V42/CTAC99/Mack}.

\bibitem{MacKenzie03}
T.~MacKenzie and A.~J. Roberts.
\newblock Holistic discretisation of shear dispersion in a two-dimensional
  channel.
\newblock In K.~Burrage and Roger~B. Sidje, editors, {\em Proc. of 10th
  Computational Techniques and Applications Conference CTAC-2001}, volume~44,
  pages C512--C530, March 2003.
\newblock
  \url{http://anziamj.austms.org.au/ojs/index.php/ANZIAMJ/article/view/694}.

\bibitem{MacKenzie05a}
T.~MacKenzie and A.~J. Roberts.
\newblock Accurately model the {Kuramoto--Sivashinsky} dynamics with holistic
  discretisation.
\newblock {\em SIAM J.~Applied Dynamical Systems}, 5(3):365--402, 2006.
\newblock \doi{10.1137/050627733}
  \url{http://epubs.siam.org/SIADS/volume-05/art_62773.html}.

\bibitem{MacKenzie05}
Tony MacKenzie.
\newblock {\em Create accurate numerical models of complex spatio-temporal
  dynamical systems with holistic discretisation}.
\newblock PhD thesis, University of Southern Queensland, 2005.

\bibitem{MacKenzie09a}
Tony MacKenzie and A.~J. Roberts.
\newblock Computer algebra derives the slow manifold of macroscale holistic
  discretisations in two dimensions.
\newblock Technical report, The University of Adelaide, 2009.
\newblock \url{http://hdl.handle.net/2440/49292}.

\bibitem{mccormick92}
S.F. Mc{C}ormick.
\newblock {\em Multilevel projection methods for partial differential
  equations}.
\newblock Society for Industrial and Applied Mathematics ({SIAM}),
  Philadelphia, 1992.

\bibitem{npl61}
{National Physical Laboratory}.
\newblock {\em Modern Computing Methods}, volume~16 of {\em Notes on Applied
  Science}.
\newblock Her Majesty's Stationery Office, London, 2nd edition, 1961.

\bibitem{Pavliotis07}
G.~A. Pavliotis and A.~M. Stuart.
\newblock {\em Multiscale methods: averaging and homogenization}, volume~53 of
  {\em Texts in Applied Mathematics}.
\newblock Springer, 2008.

\bibitem{Press92}
W.~H. Press, S.~A. Teukolsky, W.~T. Vetterling, and B.~P. Flannery.
\newblock {\em Numerical recipes in FORTRAN. The art of scientific computing}.
\newblock CUP, 2nd edition, 1992.
\newblock \url{http://www.library.cornell.edu/nr/cbookfpdf.html}.

\bibitem{Roberts00a}
A.~J. Roberts.
\newblock A holistic finite difference approach models linear dynamics
  consistently.
\newblock {\em Mathematics of Computation}, 72:247--262, 2002.
\newblock \url{http://www.ams.org/mcom/2003-72-241/S0025-5718-02-01448-5}.

\bibitem{Roberts04d}
A.~J. Roberts and I.~G. Kevrekidis.
\newblock Higher order accuracy in the gap-tooth scheme for large-scale
  dynamics using microscopic simulators.
\newblock In Rob May and A.~J. Roberts, editors, {\em Proc. of 12th
  Computational Techniques and Applications Conference CTAC-2004}, volume~46 of
  {\em ANZIAM~J.}, pages C637--C657, July 2005.
\newblock \protect \url{http://anziamj.austms.org.au/V46/CTAC2004/Robe} [July
  20, 2005].

\bibitem{Roberts06d}
A.~J. Roberts and I.~G. Kevrekidis.
\newblock General tooth boundary conditions for equation free modelling.
\newblock {\em SIAM J.~Scientific Computing}, 29(4):1495--1510, 2007.
\newblock \doi{10.1137/060654554}.

\bibitem{Roberts96a}
A.J. Roberts.
\newblock Low-dimensional modelling of dynamics via computer algebra.
\newblock {\em Comput. Phys. Comm.}, 100:215--230, 1997.

\bibitem{Roberts01b}
A.J. Roberts.
\newblock Derive boundary conditions for holistic discretisations of {Burgers'}
  equation.
\newblock Technical report, [\url{http://arXiv.org/abs/math.NA/0106224}], 2001.

\bibitem{Roberts98a}
A.J. Roberts.
\newblock Holistic discretisation ensures fidelity to {Burgers'} equation.
\newblock {\em Applied Numerical Mathematics}, 37:371--396, 2001.

\bibitem{Samaey03b}
G.~Samaey, I.~G. Kevrekidis, and D.~Roose.
\newblock The gap-tooth scheme for homogenization problems.
\newblock {\em SIAM Multiscale Modeling and Simulation}, 4:278--306, 2005.
\newblock \doi{10.1137/030602046}.

\bibitem{vorst03}
H.A. van~der Vorst.
\newblock {\em Iterative {Krylov} methods for large linear systems}, volume~13
  of {\em Cambridge monographs on applied and computational mathematics}.
\newblock Cambridge University Press, 2003.

\bibitem{Vanderbauwhede89}
A.~Vanderbauwhede.
\newblock Centre manifolds.
\newblock {\em Dynamics Reported}, pages 89--169, 1989.

\bibitem{Wolfram96}
S.~Wolfram.
\newblock {\em The Mathematica Book, 3rd ed.}
\newblock Wolfram Media/Cambridge University Press, 1996.

\end{thebibliography}

\end{document}